\documentclass[runningheads]{llncs}
\usepackage{etoolbox}
\usepackage{array,multirow}
\usepackage{graphicx}
\usepackage{amsmath}
\usepackage{amsfonts}
\usepackage{amsthm}
\usepackage{amssymb}
\usepackage{mathabx}
\usepackage{xcolor}
\usepackage{bm}
\usepackage{bbm}
\usepackage{subcaption}
\usepackage{wasysym}
\usepackage{hyperref}
\usepackage[linesnumbered,ruled]{algorithm2e}
\usepackage{adjustbox}
\usepackage{mathbbol}
\usepackage{authblk}
\usepackage{amsthm}
\usepackage{pstricks}
\usepackage{pst-eucl}
\usepackage{pst-solides3d}
\usepackage{pst-3dplot}
\usepackage{textcomp}
\usepackage{tikz}

\newcommand{\ste}[1]{\textcolor{blue}{#1}}

\def \D{\mathcal{D}_n}

\begin{document}
%\newtheorem{theorem}{Theorem}
%\newtheorem{definition}[theorem]{Definition}
%\newtheorem{lemma}[theorem]{Lemma}
%\theoremstyle{remark}
%\newtheorem{remark}[theorem]{Remark}
%\newtheorem{example}[theorem]{Example}
%
%
% \title{Three New Equivalent Encoding for Rhythmic Tiling Problems in Music}

\title{A SAT Encoding to Compute Aperiodic Tiling Rhythmic Canons}

\author{Gennaro Auricchio\inst{1}\texorpdfstring{\href{https://orcid.org/0000-0002-4285-8887}{\protect\includegraphics[scale=0.75]{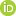}}{$\;$  }} \and
Luca Ferrarini\inst{1,2,3}\texorpdfstring{\href{https://orcid.org/0000-0002-8903-2871}{\protect\includegraphics[scale=0.75]{orcid.png}}{$\;$  }} \and
Stefano Gualandi\inst{1}\texorpdfstring{\href{https://orcid.org/0000-0002-2111-3528}{\protect\includegraphics[scale=0.75]{orcid.png}}{$\;$  }} \and \\
Greta Lanzarotto\inst{1,2,3,4}\texorpdfstring{\href{https://orcid.org/0000-0001-9743-0573}{\protect\includegraphics[scale=0.75]{orcid.png}}{$\;$  }}\and
Ludovico Pernazza\inst{1}%\orcidID{0000-0000-0000-0000}}
}
\authorrunning{G. Auricchio et al.}

\institute{University of Pavia, Department of Mathematics\and University of Milano-Bicocca, Department of Mathematics and its Applications\and INdAM \and
University of Strasbourg, IRMA}
\maketitle              % typeset the header of the contribution
%

         % typeset the header of the contribution
%
\begin{abstract}
% The Rhythmic Tiling Problem consists in finding all the possible aperiodic complements of a given rhythm. In this paper, we propose three encoding of this problem. The first one is an improved version of the already existing ILP model implemented through Gurobi. The second one is a Pseudo-Boolean formulation implemented through MiniZinc. The third one is a complete SAT formulation implemented through PySAT. We compare our models with the state of the art and we validate our encoding by testing them on several rhythms. Finally, we provide the yet unknown complete list of tiling complements of $A=\{0,20,40,45,65,85\}\subset\mathbb{Z}_{180}$.

In Mathematical Music theory, the Aperiodic Tiling Complements Problem 
% Aperiodic Tiling Rhythmic Problem 
consists in finding all the possible aperiodic complements of a given rhythm $A$.
The complexity of this problem depends on the size of the period $n$ of the canon and on the cardinality of the given rhythm $A$.
The current state-of-the-art algorithms can solve instances with $n$ smaller than $180$.
In this paper we propose an ILP formulation and a SAT Encoding to solve this mathemusical problem, and we use the Maplesat solver to enumerate all the aperiodic complements.
We validate our SAT Encoding using several different periods and rhythms
% .
% With our approach,
and we compute for the first time the complete list of aperiodic tiling complements of standard Vuza rhythms for canons of period $n=\{180,420,900\}$. 


\end{abstract}

\section{Introduction}

Mathematical Music Theory is the study of Music from a mathematical point of view. Many connections have been discovered, some of which have a long tradition, but they seem to be still offering new problems and ideas to researchers,
%{\bf TODO: manca una prima frase per dire di cosa si occupa  "the mathematical music theory" citata all'inizio dell'abstract.}
%
%{\bf TODO: io spiegherei meglio l'interplay tra matematica e musica come nell'esempio 1 del lavoro di Andreatta}
%
whether they be music composers or computer scientists.
%are interested in formalizing music problems using modern computational models.
%
The first attempt to produce music through a computational model dates back to $1957$, when the authors %{\bf TODO: dire cosa hanno fatto} 
composed a string quartet, also known as the Illiac Suite, through 
% the generate-and-test problem-solving approach
random number generators and Markov chains \cite{hiller1957musical}. 
Since then, a plethora of other works have explored how computer science and music can interact: to compose music \cite{wiggins1989representing,zimmermann2001modelling}, to analyse existing compositions and melodies \cite{chemillier2001two,courtot1990constraint,ebciouglu1988expert}, or even to represent human gestures of the music performer \cite{radicioni2007constraint}. 
In particular, Constraint Programming has been used to model harmony, counterpoint and other aspects of music (e.g., see \cite{anders2011constraint}), to compose music of various genres as described in the book \cite{anders2018compositions}, or to impose musical harmonization constraints in \cite{pachet2001musical}.
%
%\genna{Another field of research makes use of Artificial Intelligence to translate music scores into real music. 
%In \cite{de2002ai}, for example, the authors used a case-based reasoning system to recreate the expressiveness of the performer, while in \cite{wang2019performancenet} a deep convolutional model was proposed.}
%Remarkably, {\bf TODO: non basta citare, bisogna dire cosa fanno con una frase, altrimenti non serve a nulla, bisogna esser un po' più specifici (NON MI SEMBRA PROPRIO CHE PARLINO DI MACHINE LEARNING!} several studies were pursued through Machine Learning based techniques 
% \genna{In}  \genna{the authors applied AI techniques to recreate the expressiveness of musician composers}, \cite{numao2002constructive,papadopoulos1999ai}. \genna{Recently, it has also been proposed a deep convolutional model for
% translating music scores into music  \cite{wang2019performancenet}. } 

%{\bf TODO: manca completamente una introduzione di background che permetta al lettore, completamente digiuno di matemusica, di capire di cosa state parlando.}
In this paper, we deal with Tiling Rhythmic Canons, that are purely rhythmic contrapuntal compositions. % in which there are no superimposition of the voices and there are no ``holes'' \cite{Messiaen}.
%
%This problem can be expressed through algebraic tools. 
For a fixed period $n$, a tiling rhythmic canon is a couple of sets $A,B\subset\{0,1,2,\dots,n-1\}$ such that at every instant there is exactly one voice playing; $A$ defines the sequence of beats played by every voice, $B$ the instants at which voices start to play.
% 
% rhythm $A$ is a set of natural numbers, taken modulo $n$. The set $A$ represents the instants at which an instrument plays.  
% In this framework, a complement of $A$ is another set of instants, which we denote $B$, that describes the instants at which the rhythm $A$ is reproduced by a musical instrument. 
% 
%In this framework, if
% the couple $(A,B)$ is said to be a \emph{tiling rhythmic canon}. 
%Algebraically this is equivalent to require that each translation of $A$ by $b\in B$ is disjoint from the other possible translations and the union of all the translation gives the whole set of instants. This can be expressed through the identity $A\oplus B =\mathbb{Z}_n$.
If one of the sets, say $A$, is given, it is well-known that the problem 
of finding a \emph{complement} $B$ has in general no unique solution. It is very easy to find tiling canons in which at least one of the set is \emph{periodic}, i.e. it is built repeating a shorter rhythm.
From a mathematical point of view, the most interesting canons are therefore those in which both sets are \emph{aperiodic} (the problem can be equivalently rephrased as a research of tessellations of a special kind). 
%A complement is said aperiodic if it has no repeated inner structures, that is, it is obtained as a repetition of a shorter rhythm.
%
%
%More formally, the Aperiodic Tiling Complement Problem we consider is defined as follows.
%
%{\bf TODO: questa definizione è incomprensibile per chi non ha studiato algebra. Forse non la metterei nell'introduzione.}
%Given a period $n\in\mathbb{N}$ and a rhythm $A\subset \mathbb{Z}_n$, we want to determine all the subsets $B\subset \mathbb{Z}_n$ such that $A\oplus B =\mathbb{Z}_n$. 
%
%Here, $\mathbb{Z}_n$ is the cyclic group of order $n$, that is, {\bf TODO: ci vuole una spiega for dummies mirata a questo lavoro}.
%
%From a musical point of view, $A$ is the ``inner voice'' of a canon, while $B$ represents the ``outer voice''. 
%
%Given a set $A$, there might exist more than one complement, and among all the possible solutions, the aperiodic {\bf non avete detto cosa sono le aperiodiche} canons are of major interest from a mathematical point of view.  
%
Enumerating all aperiodic tiling canons has to face two main hurdles: on one side, the problem lacks the structure of other algebraic ones, such as ring or group theory; on the other side, the combinatorial size of the domain becomes enormous very soon.
%Due to the complexity of the problem {\bf non hai spiegato perchè è complesso},
Starting from the first works in the 1940s, research has gradually shed some light on parts of the problem from a theoretical point of view, and several heuristics and algorithms that allow to compute tiling complements have been introduced, but a complete solution appears to still be out of reach.

\paragraph{Contributions.}
The main contributions of this paper are the Integer Linear Programming (ILP) model and the SAT Encoding to solve the Aperiodic Tiling Complements Problem presented in Section 3.
Using a modern SAT solver we are able to compute the complete list of aperiodic tiling complements of a class of Vuza rhythms for periods $n = \{ 180, 420, 900\}$.
%
%This result is enabling new studies of the mathematical structure of Vuza rhythms, keeping alive the dynamic movements between musical problems and mathematical statements.

\paragraph{Outline.} 
The outline of the paper is as follows. 
Section \ref{sec:basic_notions} reviews the main notions on Tiling Rhythmic Canons and defines formally the problem we tackle.
% of finding aperiodic complements of a given rhythm. 
%
In Section \ref{sec:our_contribution}, we introduce an ILP model and a SAT Encoding of the Aperiodic Tiling Complements Problem expressing the tiling and the aperiodicity constraints in terms of Boolean variables.
Finally, in Section \ref{sec:final}, we include our computational results to compare the efficiency of the aforementioned ILP model and SAT Encoding with the current state-of-the-art algorithms.
%as explained in \cite{andreatta2011constructing}.

\section{The Aperiodic Tiling Complements Problem}
\label{sec:basic_notions}
We begin fixing some notation and giving the main definitions.
In the following, we conventionally denote the cyclic group of remainder classes modulo $n$ by $\mathbb{Z}_n$ and its elements with the integers $\{0, 1, \dots, n - 1 \}$, i.e. identifying each class with its least non-negative member.
\begin{definition}\label{directsum}
    Let $A, B \subset \mathbb{Z}_n$. Let us define the application
\[\sigma:A \times B \rightarrow \mathbb{Z}_n, (a, b) \mapsto a + b.\]
We set $A + B: = \mbox{Im}(\sigma)$; if $\sigma$ is bijective we say that $A$ and $B$ {\bf are in direct sum}, and we write
\[A \oplus B: =  \mbox{Im}(\sigma).\]
If $\mathbb{Z}_n = A\oplus B$, we call $(A, B)$ a {\bf tiling rhythmic canon} of {\bf period $n$}; $A$ is called the {\bf inner voice} and $B$ the {\bf outer voice} of the canon.
\end{definition}
%\begin{definition}\label{df:TRC}
		%A {\bf tiling rhythmic canon} $(A,B)$ with {\bf period} $n$ is a factorization of the cyclic group $\mathbb{Z}_{n}$ with two subsets $A$ and $B$:
	%\[
	 %   A\oplus B=\mathbb{Z}_{n},
	%\]
	%where $\oplus$ denotes the group operation $A\oplus B := \{([a]_n+[b]_n) \mid [a]_n \in A, [b]_n \in B\}$.
	%
	%We say that $A$ tiles with $B$, and vice versa.
%\end{definition}

%Let us pin down what is meant by identical canons. T
\begin{remark} It is easy to see that the tiling property is invariant under translations, i.e. if $A$ is a tiling complement of some set $B$, also any translate $A + z$ of $A$ is a tiling complement of $B$ (and any translate of $B$ is a tiling complement of $A$). In fact, suppose that $A \oplus B = \mathbb{Z}_n$; for every $k, z \in \mathbb{Z}_n$ %we have that $k - z \in \mathbb{Z}_n$ therefore 
by definition there exists one and only one pair $(a,b) \in A\times B$ such that $k - z = a + b$. Consequently, there exists one and only one pair $(a + z, b) \in (A + z)\times B$ such that $k = (a + z) + b$, that is $(A + z)\oplus B =\mathbb{Z}_n$.
In view of this, without loss of generality, we shall limit our investigation to rhythms containing 0 and consider equivalence classes under translation. 
\end{remark}

%Notice that while in mathematics $A \oplus B = B \oplus A$, from a musical perspective, there is a difference between canon $(A,B)$ and canon $(B,A)$: the listener recognizes the inner rhythm (the first element of the couple) performed by a certain voice and repeated by further instruments starting at the instants indicated by the outer rhythm (the second element).
%
%For a pictorial example see Figures~\ref{fig:A} and \ref{fig:B}.
%
%{\bf Ma perchè lo facciamo notare? Serve per capire il resto dell'articolo? me lo sono perso.}

\begin{figure}[t!]
	\centering
	\begin{minipage}[b]{0.45\linewidth}
		\centering	
					\begin{tikzpicture}[scale=0.5]
			\centering
			\draw[step=1cm] (0,0) grid (9,3);
			\fill[black] (0,20mm) -- (10mm,20mm) -- (10mm,30mm)-- (0,30mm) node[above right] {\color{black}0};
			\fill[black] (0,20mm) -- (10mm,20mm) -- (10mm,30mm)-- (0,30mm) node[above right] {\color{black}};
			\fill[black] (10mm,20mm) -- (20mm,20mm) -- (20mm,30mm)-- (10mm,30mm) node[above right] {\color{black}1};
			\draw (20mm,20mm) -- (30mm,20mm) -- (30mm,30mm)-- (20mm,30mm) node[above right] {\color{black}2};
			\fill[black] (30mm,10mm) -- (40mm,10mm) -- (40mm,20mm)-- (30mm,20mm);
			\draw(30mm,20mm) -- (40mm,20mm) -- (40mm,30mm)-- (30mm,30mm)node[above right] 	{\color{black}3};
			\draw(40mm,20mm) -- (50mm,20mm) -- (50mm,30mm)-- (40mm,30mm)node[above right] 	{\color{black}4};
			\fill[black] (40mm,10mm) -- (50mm,10mm) -- (50mm,20mm)-- (40mm,20mm);
			\fill[black] (50mm,20mm) -- (60mm,20mm) -- (60mm,30mm)-- (50mm,30mm) node[above right] {\color{black}5};
			\fill[black] (60mm,0) -- (70mm,0) -- (70mm,10mm)-- (60mm,10mm);
			\fill[black] (70mm,0) -- (80mm,0) -- (80mm,10mm)-- (70mm,10mm);
			\fill[black] (20mm,0) -- (30mm,0) -- (30mm,10mm)-- (20mm,10mm);
			\fill[black] (80mm,10mm) -- (90mm,10mm) -- (90mm,20mm)-- (80mm,20mm);
			\draw (60mm,20mm) -- (70mm,20mm) -- (70mm,30mm)-- (60mm,30mm)node[above right] {\color{black}6};
			\draw (70mm,20mm) -- (80mm,20mm) -- (80mm,30mm)-- (70mm,30mm)node[above right] {\color{black}7};
			\draw (80mm,20mm) -- (90mm,20mm) -- (90mm,30mm)-- (80mm,30mm)node[above right] {\color{black}8};
		\end{tikzpicture}
		\caption{$A=\{0, 1, 5\}$ as inner voice. }
		\label{fig:A}
	\end{minipage}
	\quad
	\begin{minipage}[b]{0.45\linewidth}
		\centering
		\begin{tikzpicture}[scale=0.5]
			\centering
			\draw[step=1cm] (0,0) grid (9,3);
			\fill[black] (0,20mm) -- (10mm,20mm) -- (10mm,30mm)-- (0,30mm) node[above right] {\color{black}0};
			%\draw[black] (10,20mm) -- (20mm,20mm) -- (20mm,30mm)-- (10,30mm) node[above right] {\color{black}1};
			\draw[black] (10mm,20mm) -- (20mm,20mm) -- (20mm,30mm)-- (10mm,30mm) node[above right] {\color{black}1};
			\draw (20mm,20mm) -- (30mm,20mm) -- (30mm,30mm)-- (20mm,30mm) node[above right] {\color{black}2};
			\fill[black] (40mm,10mm) -- (50mm,10mm) -- (50mm,20mm)-- (40mm,20mm);
			\fill(30mm,20mm) -- (40mm,20mm) -- (40mm,30mm)-- (30mm,30mm)node[above right] 	{\color{black}3};
			\draw(40mm,20mm) -- (50mm,20mm) -- (50mm,30mm)-- (40mm,30mm)node[above right] 	{\color{black}4};
			\fill(10mm,10mm) -- (20mm,10mm) -- (20mm,20mm)-- (10mm,20mm);
			\draw(50mm,20mm) -- (60mm,20mm) -- (60mm,30mm)-- (50mm,30mm) node[above right] {\color{black}5};
			\fill[black] (50mm,0) -- (60mm,0) -- (60mm,10mm)-- (50mm,10mm);
			\fill[black] (80mm,0) -- (90mm,0) -- (90mm,10mm)-- (80mm,10mm);
			\fill[black] (20mm,0) -- (30mm,0) -- (30mm,10mm)-- (20mm,10mm);
			\fill[black] (70mm,10mm) -- (80mm,10mm) -- (80mm,20mm)-- (70mm,20mm);
			\fill (60mm,20mm) -- (70mm,20mm) -- (70mm,30mm)-- (60mm,30mm)node[above right] {\color{black}6};
			\draw (70mm,20mm) -- (80mm,20mm) -- (80mm,30mm)-- (70mm,30mm)node[above right] {\color{black}7};
			\draw (80mm,20mm) -- (90mm,20mm) -- (90mm,30mm)-- (80mm,30mm)node[above right] {\color{black}8};
		\end{tikzpicture}
		\caption{$B=\{0,3,6\}$ as inner voice.}
		\label{fig:B}
	\end{minipage}
\end{figure}

\begin{example}
We consider a period $n=9$, %and the cyclic group $\mathbb{Z}_9=\{0,1, \dots, n-1\}$.
and the two rhythms $A=\{0,1,5\} \subset \mathbb{Z}_9$ and $B=\{0,3,6\} \subset \mathbb{Z}_9$ in Figure \ref{fig:A} and Figure \ref{fig:B}. 
They provide the canon $A \oplus B = \mathbb{Z}_9$, since
$\{0,1,5\}  \oplus \{0,3,6\} = \{0, 3, 6, 1, 4, 7, 5, 8, 2\}$, where the last number is obtained by $(5+6) \mod 9 = 2$. 
\end{example}

\begin{definition}\label{def:period}
    A rhythm $A \subset\mathbb{Z}_n$ is {\bf periodic (of period $z$)} if and only if there exists an element $z \in \mathbb{Z}_n$, $z\neq 0$, such that $z + A = A$. 
    In this case, $A$ is also called periodic modulo $z\in\mathbb{Z}_n$.
    A rhythm $A\subset\mathbb{Z}_n$ is {\bf aperiodic} if and only if it is not periodic.
\end{definition} 

Coming back to Example 1, it is easy to note the periodicity $z = 3$ in rhythm $B=\{0, 3, 6\}$: indeed, $3 + B = B$.
Notice that if $A$ is periodic of period $z$, $z$ must be a strict divisor of the period $n$ of the canon.

Tiling rhythmic canons can be characterised using polynomials, as follows.
	
\begin{lemma}
    \label{lm:pol_equivalence}
    Let $A$ be a rhythm in $\mathbb{Z}_n$ and let $p_A(x)$ be the {\bf characteristic polynomial} of $A$, that is, $p_A(x)=\sum_{k\in A}x^{k}$. Given $B\subset\mathbb{Z}_n$ and its characteristic polynomial $p_B(x)$, we have that
\begin{equation}\label{eq:pol_form}
p_A (x)\cdot p_B (x)\equiv \sum_{k=0}^{n-1} x^k,\quad\quad\mod (x^{n} - 1)   
\end{equation}
 if and only if $p_A (x), p_B (x)$ are polynomials with coefficients in $\{0,1\}$ and $A\oplus B = \mathbb{Z}_n$.
\end{lemma}

%\begin{definition}
    %Let $k>0, k  \in \mathbb{Z}_{n}$ be a divisor of $n$. A rhythm $A\subset \mathbb{Z}_{n}$ is {\bf periodic modulo} $k$ if and only if  $\{k\} \oplus A = A$ ({\bf dove adesso la somma viene fatta...}). A rhythm $A\subset \mathbb{Z}_{n}$ is {\bf aperiodic} if and only if it is not periodic for any $k \in \mathbb{Z}_n$ such that $k$ is a divisor of $n$.
%\end{definition}

%
%This implies the identification of an internal period, which divides the period of the canon.
%
%In this case, $3 + B = B$ and hence $B$ is periodic modulo \ste{$z=3$ (z invece di k)}.

\begin{definition}
    A tiling rhythmic canon $(A,B)$ in $\mathbb{Z}_{n}$ is a {\bf Vuza canon} if both $A$ and $B$ are aperiodic.
\end{definition}

\begin{remark}
    \label{aperiodic}
    Note that a set $A$ is periodic modulo $z$ if and only if it is periodic modulo all the non-trivial multiples of $z$ dividing $n$. 
    For this reason, when it comes to check whether $A$ is periodic or not, it suffices to check if $A$ is periodic modulo $m$ for every $m$ in the set of maximal divisors of $n$.
    We denote by $\D$ this set:% of maximal divisors of $n$:
\begin{equation*}
    \D:=\big\{n/ p \mid p \mbox{ is a prime factor of } n\big\}.
    \end{equation*}
We also denote with $k_n$ the cardinality of $\D$, so that $n=p_1^{\alpha_1}p_2^{\alpha_2}\dots p_{k_n}^{\alpha_{k_n}}$ is the unique prime factorization of $n$, where $\alpha_1,\dots,\alpha_{k_n}\in\mathbb{N^+}$ .
% {\bf where $\alpha$ is ...}
\end{remark}

%Given a fixed period $n\in\mathbb{N}$ and a rhythm $A$, a classical musical problem is to determine whether a rhythm $B$ which tiles with $A$ exists. 
%
For a complete and exhaustive discussion on tiling problems, we refer the reader to \cite{amiot2011structures}.
In this paper, we are interested in the following tiling problem.

\begin{definition}
    Given a period $n\in\mathbb{N}$ and a rhythm $A \subset \mathbb{Z}_n$, the {\bf Aperiodic Tiling Complements Problem} consists in finding all aperiodic complements $B$ i.e., subsets $B$ of $\mathbb{Z}_n$ such that $A \oplus B = \mathbb{Z}_n$.
\end{definition}
   
%{\bf TODO: manca la complessità teorica del problema, me ne avevate parlato} This problem is NP-complete, as proven in ??.
Some problems very similar to the decision of tiling (i.e., the tiling decision problem DIFF
in \cite{Matolcsi}) have been shown to be NP-complete; a strong lower bound for computational complexity of the tiling decision problem is to be expected, too. 
% We conjecture therefore that the problem is NP-complete in general.

%%%%%%%%%%%%%%%%%%%%%%%%%%%%%%%%%%%%%%%%%%%%%%%%%%%%
\section{A SAT Encoding }
\label{sec:our_contribution}

In this section, we present in parallel an ILP model and a new SAT Encoding for the Aperiodic Tiling Complements Problem that are both  used to enumerate all complements of $A$.
%(ERA: satisfiable models of the given CNF formula)
%
%To define our \ste{formulations}, we have to 
We define two sets of constraints: (i) the {\it tiling constraints} that impose the condition $A \oplus B = \mathbb{Z}_n$, and (ii) the {\it aperiodicity constraints} that impose that the canon $B$ is aperiodic.

\paragraph{Tiling constraints.}
Given the period $n$ and the rhythm $A$, let $\bm a=[a_0,\dots,a_{n-1}]^\intercal$ be its characteristic (column) vector, that is, $a_i=1$ if and only if $i \in A$. 
Using vector $\bm a$ we define the circulant matrix $T \in \{0,1\}^{n \times n}$ of rhythm $A$, that is, each column of $T$ is the circular shift of the first column, which corresponds to vector $\bm a$.
Thus, the matrix $T$ is equal to 
\[ 
T=
\begin{bmatrix}
    a_{0} & a_{n-1} & a_{n-2} & \dots  & a_{1} \\
    a_{1} & a_{0} & a_{n-1} & \dots  & a_{2} \\
    \vdots & \vdots & \vdots & \ddots & \vdots \\
    a_{n-1} & a_{n-2} & a_{n-3} & \dots  & a_{0}
\end{bmatrix}.
\]
We can use the circulant matrix $T$ to impose the tiling conditions as follows.
Let us introduce a literal $x_i$ for $i=0,\dots,n-1$, that represents the characteristic vector of the tiling rhythm $B$, that is, $x_i = 1$ if and only if $i \in B$.
Note that a literal is equivalent to a 0--1 variable in ILP terminology.
Then, the tiling condition can be written with the following linear constraint:
\begin{equation}\label{complementary}
    \sum_{i \in \{0, \dots, n-1\}} T_{ij} x_i = 1, \quad \forall j = 0, \dots, n-1.
\end{equation}
Notice that the set of linear constraints \eqref{complementary} imposes that exactly one variable (literal) in the set $\{x_{{n+i-j \mod n}}\}_{j\in A}$ is equal to one. Hence, we encode this condition as an 
{\tt Exactly-one} constraint, that is, exactly one literal can take the value one.
The {\tt Exactly-one} constraint can be expressed as the conjunction of the two constraints {\tt At-least-one} and {\tt At-most-one}, for which standard SAT encoding exist (e.g., see \cite{bailleux2003efficient,philipp2015pblib}). 
%
%The conjunction of the last two constraints are translated into the  $\phi= C_l \land C_m$, where $C_l = (x_{0} \lor x_{1} \lor \dots \lor x_{n-1})$ and $C_m=\bigwedge_{k\neq l}(\lnot x_{k} \lor \lnot x_{l})$ {\bf TODO: contrallare bene gli indici}
%
Hence, the tiling constraints \eqref{complementary} are encoded with the following set of clauses depending on $i=0, \dots, n-1$:
\begin{equation}\label{sat:compl}
\bigvee_{j \in A}\left(x_{n-(j-i) \mod n}\right) \bigwedge_{k,l \in A,k \neq l}\left(\lnot x_{{n-(k-i) \mod n}} \lor \lnot x_{{n-(l-i) \mod n}}\right).
\end{equation}

\paragraph{Aperiodicity constraints.}
%{\bf TODO IMPORTANT: devo rivedere bene gli indici in base a come sono definiti in $\D$. Sono partito dal codice e ricostruito il testo, non è detto che sia corretto.}
%To impose the aperiodicity constraint, we exploit the set of maximal divisor $\D$ adapting the linear constraints introduced in \cite{DBLP:journals/corr/abs-2107-04108}, which we review next.
%
In view of Definition \ref{def:period},
% to impose that $B$ is not periodic of period $d$ it is sufficient to check whether $(d+b) \mod n \in B$ for every $b\inB$:
%
%Given a fixed $d$, whenever 
if there exists a $b \in B$ such that $(d + b) \mod n \neq b$, then the canon $B$ is not periodic modulo $d$.
Notice that by Remark \ref{aperiodic} we need to check this condition only for the values of $d \in \D$.

We formulate the aperiodicity constraints introducing auxiliary variables $y_{d,i},z_{d,i},u_{d,i} \in \{0,1\}$ for every prime divisor $d \in \D$ and for every integer $i = 0,\dots,d-1$.
We set
% \small
% \begin{align}
% \label{implications} &u_{d,i} = 1 \; \Leftrightarrow \;
%     \left(\sum_{k=0}^{n/d-1} x_{i+kd} = \frac{n}{d}\right) \vee \left(\sum_{k=0}^{n/d-1} x_{i+kd} = 0\right), &\forall d \in \D, i=0,\dots,d-1,\\
% \label{sumdivisor}     &\sum_{i=0}^{d-1} u_{d,i} \leq d - 1 &\forall d \in \D.
% \end{align}
% \normalsize
\begin{equation}
    \label{implications} u_{d,i} = 1 \; \Leftrightarrow \;
    \left(\sum_{k=0}^{n/d-1} x_{i+kd} = \frac{n}{d}\right) \vee \left(\sum_{k=0}^{n/d-1} x_{i+kd} = 0\right),
\end{equation}
for all  $d \in \D$, $i=0,\dots,d-1$, with the condition
\begin{equation}
\label{sumdivisor}   
\sum_{i=0}^{d-1} u_{d,i} \leq d-1, \quad \forall d \in \D.
\end{equation}

Similarly to \cite{auricchio2021integer}, the constraints \eqref{implications} can be linearized using standard reformulation techniques as follows:

\begin{align}
\label{y:1}     & 0 \leq \sum_{k=0}^{n/d} x_{i+kd} - \frac{n}{d}y_{d,i}\leq \frac{n}{d} - 1 & \forall d \in \D,\;\; i=0,\dots,d-1, \\
\label{z:-1}    & 0 \leq  \sum_{k=0}^{n/d} (1-x_{i+kd}) - \frac{n}{d}z_{d,i} \leq \frac{n}{d} - 1 & \forall d \in \D, \; \; i=0,\dots,d-1,\\
\label{U}       & y_{d,i} + z_{d,i} = u_{d,i} & \forall d \in \D, \;\; i=0,\dots,d-1. 
\end{align}

\noindent Notice that when $u_{d,i}=1$ exactly one of the two incompatible alternatives in the right hand side of \eqref{implications} is true, %are equivalent to a single equality constraint, 
while whenever $u_{d,i}=0$ the two constraints are false. 
%
%\genna{The constraint \eqref{U} imposes that if one of the two variables $y_{d,i}$ and $z_{d,i}$ is equal to one, the corresponding auxiliary variable $u_{d,i}$ is set to one.}
%
Correspondingly, the constraint \eqref{U} imposes that the variables $y_{d,i}$ and $z_{d,i}$ cannot be equal to $1$ at the same time.
On the other hand, constraint \eqref{sumdivisor} imposes that at least one of the auxiliary variables $u_{d,i}$ be equal to zero.
%Let $n = p_{1}^{\alpha_1}p_{2}^{\alpha_2}\dots p_{K_n}^{\alpha_{K_n}}$ be the prime decomposition of $n$. In the ILP model, to impose the aperiodicity constraints, we introduced a family of Boolean variables $\Uij$, where 
% $m_{j}\in D_n$
%$j=1,..,K_n$ and $i = 0, \dots m_{j}-1$.
% we introduce the auxiliary Boolean variables . 
%Each variable $\Uij$ is equal to $1$, if and only if all the variables $\{x_{i+km_j}\}_{k=0,\dots,p_{j}-1}$ are all either equal to $1$ or all equal to $0$. 
% \forall i= 0,..,m_{j}$. w
%The aperiodicity of the rhythm is then imposed by constraining the number of non-null variables among the $\Uij$'s for any fixed $j$, as in \eqref{m1:c7}. %In particular, given $m\in \D$, we must have, at most, $\frac{n_B n}{m}$, where $n_B$ is the expected cardinality of the generic complement $B$.
%
%{\bf Rimetti la spiega di Luca, ma sistemata in base a quanto fatto perl'ILP.}
%

Next, we encode the previous conditions as a SAT formula.
%
%\genna{To enforce the constraint that the set of variables must be of the same type we apply the {\tt Same} function}
%
To encode the if and only if clause, we make use of the logical equivalence between $C_1 \Leftrightarrow C_2$ and $(\lnot C_1 \lor C_2) \land (C_1 \lor \lnot C_2)$.
The clause $C_1$ is given directly by the literal $u_{d,i}$.
The clause $C_2$, expressing the right hand side of \eqref{implications}, i.e. the constraint that the variables must be either all true or all false, can be written as
%
%Thus, \eqref{implications} can be written as the formula 
\[
C_2 = \left(\bigwedge_{k=0}^{n/d} x_{i+kd}\right) \vee \left(\bigwedge_{k=0}^{n/d} \bar{x}_{i+kd}\right), \quad \forall d \in \D.
\]
Then, the linear constraint \eqref{sumdivisor} can be stated as the SAT formula:
\[
  \lnot \left(u_{d,0} \land u_{d,1} \land \dots \land u_{d,(d-1)}\right) = \bigvee_{l=0}^{d-1} \bar{u}_{d,l}, \quad 
  \forall d \in \D.
\]
Finally, we express the aperiodicity constraints using
\begin{equation}\label{sat:apreriodic}
    \bigwedge\limits_{i = 0}^{d-1}
     \left[\left( \lnot C_2 \lor u_{d,i} \right)\land
    \left( C_2 \lor \bar{u}_{d,i} \right) \right]
    \land 
    \bigvee_{l=0}^{d-1} \bar{u}_{d,l},\,
    \forall d \in \D.
\end{equation}
Note that joining \eqref{complementary}, \eqref{y:1}--\eqref{U} with a constant objective function gives a complete ILP model, which can be solved with a modern ILP solver such as Gurobi to enumerate all possible solutions.
At the same time, joining \eqref{sat:compl} and \eqref{sat:apreriodic} into a unique CNF formula, we get our complete SAT Encoding of the Aperiodic Tiling Complements Problem.
(see Section 4 for computational results).

%%%%%%%%%%%%%%%%%%%%%%%%%%%%%%%%%%%%%%%%%%%%%%%%%%%%%%%%%%%    
\subsection{Existing solution approaches}

%Despite its simple formulation, 
For the computation of all the aperiodic tiling complements of a given rhythm
%is a challenging task, and we describe below 
the two most successful approaches already known are the \emph{Fill-Out Procedure} \cite{kolountzakis2009algorithms} and the {\it Cutting Sequential Algorithm} \cite{auricchio2021integer}.
%

%Notice that since any translation of a solution (canon) is still a solution \cite{Vuza}, 
%Therefore, the output of the algorithms are usually normalized (i.e., $0$ belongs to every computed canons {\bf questa non la capisco}) and all the translations of the solutions are identified in postprocessing but only one is retained.

\paragraph{The Fill-Out Procedure.} %{\bf (Da snellire di un paio di frasi)}
%
%The \emph{Fill-Out Procedure} is a heuristic algorithm introduced in \cite{kolountzakis2009algorithms} to find a complete tiling of $\mathbb{Z}_{144}$. 
%
%The key idea behind this algorithm is the following. 
%
%Given a rhythm $A\subset\mathbb{Z}_n$ such that $0\in A$, the algorithm sets $P=\{0\}$ and starts the search for possible expansions of the set $P$. 
%
%The expansion is accomplished by adding, whenever it is possible, an element $\alpha\in\mathbb{Z}_n$ to $P$. 
%
%To decide which element of $\mathbb{Z}_n$ to add, a ranking function $r(x, P)$ is introduced. 
%
%The function $r$ counts all the possible ways that $x$ can be covered through a translation of $A$ and, once every element of $\mathbb{Z}_m\backslash (A\oplus P)$ has been ranked, the algorithm tries to add the element with the lowest rank. 
%
%This defines a new set, $\Tilde{P}\supset P$, which is again expanded until either it can no longer be expanded or the set becomes a tiling complement. 
%
%The search ends when all the possibilities have been explored. 
%The algorithm finds also periodic solutions. 
%
%These solutions, alongside the translations, are removed post-processing.

The \emph{Fill-Out Procedure} is the heuristic algorithm introduced in \cite{kolountzakis2009algorithms}. 
The key idea behind this algorithm is the following: 
given a rhythm $A\subset\mathbb{Z}_n$ such that $0\in A$, the algorithm sets $P=\{0\}$ and starts the search for possible expansions of the set $P$. 
The expansion is accomplished by adding an element $\alpha\in\mathbb{Z}_n$ to $P$ according to the reverse order induced by a ranking function $r(x, P)$, which counts all the possible ways in which $x$ can be covered through a translation of $A$.
% Once every element of $\mathbb{Z}_m\backslash (A\oplus P)$ has been ranked, the algorithm tries to add the element with the lowest rank. 
%
This defines a new set, $\Tilde{P}\supset P$, which is again expanded until either it can no longer be expanded or the set becomes a tiling complement. 
The search ends when all the possibilities have been explored. 
The algorithm finds also periodic solutions that must removed in
post-processing, as well as multiple translations of the same rhythm.

\paragraph{The Cutting Sequential Algorithm (CSA).}
In \cite{auricchio2021integer}, the authors formulate the Aperiodic Tiling Complements Problem using an Integer Linear Programming (ILP) model that is based on the polynomial characterization of tiling canons.
The ILP model uses auxiliary 0--1 variables to encode the product $p_A (x) \cdot p_B(x)$ which characterizes tiling canons.
The aperiodicity constraint is formulated analogously to what done above.
The objective function is equal to a constant and does not impact the solutions found by the model.
The ILP model is used within a sequential cutting algorithm that adds a no-good constraint every time a new canon $B$ is found to prevent finding solutions twice.
In addition, the sequential algorithm sets a new no-good constraints for every translation of $B$; hence, in contrast to the \emph{Fill-Out Procedure}, the \emph{CSA Algorithm} does not need any post-processing.

%%%%%%%%%%%%%%%%%%%%%%%%%%%%%%%%%%%%%%%%%%%%%%%%%%%%%%%%%%
\section{Computational Results}
\label{sec:final}
First, we compare the results obtained using our ILP model and SAT Encoding with the runtimes of the \emph{Fill-Out Procedure} and of the \emph{CSA Algorithm}.
We use the canons with periods 72, 108, 120, 144 and 168 that have been completely enumerated by Vuza \cite{Vuza}, Fripertinger \cite{fripertinger2005remarks}, Amiot \cite{amiot2009new}, Kolountzakis and Matolcsi \cite{kolountzakis2009algorithms}. 
Table \ref{tab1} shows clearly that the two new approaches outperform the state-of-the-art, and in particular, that SAT provides the best solution approach. 
We then choose some periods $n$ with more complex prime factorizations, such as $n = p^2q^2r=180$, $n = p^2qrs=420$, and $n = p^2q^2r^2=900$. 
To find aperiodic rhythms $A$, we apply Vuza's construction \cite{Vuza} with different choices of parameters $p_1$, $p_2$, $n_1$, $n_2$, $n_3$. Thus, having $n$ and $A$ as inputs, we search for all the possible aperiodic complements and then we filter out the solutions under translation.
Since the post-processing is based on sorting canons, it requires a comparatively small amount of time.
We report the results in Table \ref{tab2}: the solution approach based on the SAT Encoding is the clear winner (from a Music theory perspective, it is also noteworthy that this is the first time that all the tiling complements, whose number is reported in the last column of the two tables, of the studied rhythms are computed).
% compare our models with the $FOP$ and with the $CSA$. \genna{referenza per minizinc} \cite{minizinc2007}, \genna{PySAT} \cite{pysat2018},

\paragraph{Implementation Details.}
We have implemented in Python the ILP model and in PySat \cite{imms-sat18} the SAT Encoding discussed in Section 3. We use Gurobi 9.1.1 as ILP solver and Maplesat \cite{phdthesis} as SAT solver.
The experiments are run on a Dell Workstation with a Intel Xeon W-2155 CPU with 10 physical cores at 3.3GHz and 32 GB of RAM.
In case of acceptance, we will release the source code and the instances on GitHub.
%The source code and the dataset is freely available on GitHub at \url{https://github.com/stegua/total-matching}
\subsubsection{Conclusions and Future Work.}
It is thinkable to devise an algorithm that, for a given $n$, finds all the pairs $(A, B)$ that give rise to a Vuza canon of period $n$.
This could provide in-depth information on the structure of Vuza canons.

\subsubsection{Acknowledgements.}
This research was partially supported by: Italian Ministry of Education, University and Research (MIUR), Dipartimenti di Eccellenza Program (2018--2022) - Dept. of Mathematics ``F. Casorati'', University of Pavia; Dept. of Mathematics and its Applications, University of Milano-Bicocca; National Institute of High Mathematics (INdAM) ``F. Severi''; Institute for Advanced Mathematical Research (IRMA), University of Strasbourg.
\begin{table}
%\tiny
\caption{Aperiodic tiling complements for periods $n\in\{72,108,120,144,168\}$.}
\label{tab1}
\vspace{.2cm}
\centering
\begin{adjustbox}{width=0.9\textwidth}
\begin{tabular}{|c|c|c|c|c|c|c|r|r|r|r|r|}
\hline
\multirow{2}{*}{$n$} &\multirow{2}{*}{$\D$}&\multirow{2}{*}{$p_1$}&\multirow{2}{*}{$n_1$}&\multirow{2}{*}{$p_2$}&\multirow{2}{*}{$n_2$}&\multirow{2}{*}{$n_3$}&\multicolumn{4}{c|}{runtimes (s)}& \multirow{2}{*}{$\# B$}\\
\cline{8-11}
&&&&&&& \emph{FOP}& \emph{CSA}& \emph{SAT} &  \emph{ILP} & \\
\hline
\hline
72&$\{24,36\}$&2&2&3&3&2& 1.59 &0.10 &$ < 0.01$  & 0.03 &6\\
\hline
\hline
108&$\{36,54\}$&2&2&3&3&3& 896.06 &7.84 &0.09 & 0.19  & 252\\
\hline
\hline
\multirow{2}{*}{120}&\multirow{2}{*}{$\{24,40,60\}$}&2&2&5&3&2&24.16&0.27&0.02& 0.04  &18\\
&&2&2&3&5&2&10.92&0.14&0.01& 0.04  &20\\
\hline
\hline
\multirow{4}{*}{144}&\multirow{4}{*}{$\{48,72\}$}&4&2&3&3&2&82.53&2.93&0.02& 0.11  &36\\
&&2&2&3&3&4&$> 10800.00$ &$> 10800.00$&11.04& 46.96 &8640\\
&&2&2&3&3&4&7.13&0.10&$< 0.01$& 0.05 &6\\
&&2&4&3&3&2&80.04&0.94&0.02& 0.08 &60\\
\hline
\hline
\multirow{2}{*}{168}&\multirow{2}{*}{$\{24,56,84\}$}&2&2&7&3&2&461.53&17.61&0.04& 0.20 &54\\
&&2&2&3&7&2&46.11&0.91&0.02& 0.07 &42\\
\hline
\end{tabular}
\end{adjustbox}
\vspace{1cm}
\end{table}
\begin{table}
\caption{Aperiodic tiling complements for periods $n\in\{180,420,900\}$.}
\label{tab2}
\vspace{.2cm}
\centering
\begin{adjustbox}{width=0.82\textwidth}
\begin{tabular}{|c@{ }|c@{ }|c@{ }|c@{ }|c@{ }|c@{ }|c@{ }|r@{ }|r@{ }|r@{ }|}
\hline
\multirow{2}{*}{$n$} &\multirow{2}{*}{$\D$}&\multirow{2}{*}{$p_1$}&\multirow{2}{*}{$n_1$}&\multirow{2}{*}{$p_2$}&\multirow{2}{*}{$n_2$}&\multirow{2}{*}{$n_3$}&\multicolumn{2}{c|}{runtimes (s)}& \multirow{2}{*}{$\# B$}\\
\cline{8-9}
&&&&&&& \emph{SAT} &  \emph{ILP} & \\
\hline
% $n$ &$D$&$A$ & $\ilp B$\\
% \hline
% \hline
% 72&$\{24,36\}$&$\{0, 8, 16, 18, 26, 34\}$&6\\
% \hline
% \hline
% 108&$\{36,54\}$&$\{ 0, 12, 24, 27, 39, 51\}$&252\\
% \hline
% \hline
% 120&$\{24,40,60\}$&$\{0, 8, 16, 30, 38, 46\}$&18\\
% 120&&$\{0, 8, 16, 24, 30, 32, 38, 46, 54, 62\}$&20\\
% \hline
% \hline
% 144&$\{48,72\}$&$\{{0, 16, 18, 32, 34, 50}\}$&36\\
% 144&&$\{0, 16, 32, 36, 52, 68\}$&8640\\
% 144&&$\{0, 9, 16, 25, 32, 36, 41, 45, 52, 61, 68, 77\}$&6\\
% 144&&$\{0, 16, 18, 32, 34, 36, 50, 52, 54, 68, 70, 86\}$&60\\
% \hline
% \hline
% 168&$\{24,56,84\}$&$\{0, 8, 16, 42, 50, 58\}$&54\\
% 168&&$\{0, 8, 16, 24, 32, 40, 42, 48, 50, 58, 66, 74, 82, 90\}$&42\\
% \hline
 \hline
\multirow{5}{*}{180}&\multirow{5}{*}{$\{36,60,90\}$}&2&2&5&3&3& 2.57 & 5.62 &2052\\
\cline{3-10}
&&3&3&5&2&2  &0.07 & 0.14 &96\\
\cline{3-10}
&&2&2&3&5&3   & 1.25 & 2.23  &1800\\
\cline{3-10}
&&2&5&3&3&2 & 0.05&  0.16 & 120\\
\cline{3-10}
&&2&2&3&3&5& 8079.07 & $> 10800.00$  & 281232\\
\hline
\hline
\multirow{12}{*}{420} &\multirow{12}{*}{$\{60,84,140,210\}$}&7&5&3&2&2 &2.13 & 3.57  &720 \\
\cline{3-10}
&&5&7&3&2&2 & 1.52 & 4.08  &672 \\
\cline{3-10}
&&7&5&2&3&2& 7.73 & 16.11 & 3120 \\
\cline{3-10}
&&5&7&2&3&2 & 1.63 &  4.18  & 1008 \\
\cline{3-10}
&&7&3&5&2&2 & 4.76 & 7.45  & 864 \\
\cline{3-10}
&&3&7&5&2&2 & 12.78 & 32.19  & 6720 \\
\cline{3-10}
&&7&3&2&5&2& 107.83 & 1186.21 & 33480 \\
\cline{3-10}
&&3&7&2&5&2&0.73 &  2.36 & 840\\
\cline{3-10}
& &7&2&5&3&2 & 11.14 & 21.19 & 1872 \\
\cline{3-10}
&&2&7&5&3&2& 17.31&  52.90 & 10080 \\
\cline{3-10}
&&7&2&3&5&2& 89.97 & 691.56  & 22320 \\
\cline{3-10}
&&2&7&3&5&2 & 1.17 & 4.13 & 1120 \\
\hline
\hline
\multirow{5}{*}{900}&\multirow{5}{*}{$\{180,300,450\}$}&2&25&3&3&2 & 43.60 & 110.65 & 15600 \\
\cline{3-10}
&&5&10&3&3&2 & 107.36 & 741.79 & 15840 \\
\cline{3-10}
& &2&9&5&5&2   & 958.58 & $> 10800.00$ & 118080 \\
\cline{3-10}
& & 6&3&5&5&2 &5559.76 &$> 10800.00$ &123840 \\
\cline{3-10}
&&3  & 6&5&5&2&486.39  & 8290.35 & 62160\\ 
\hline
\end{tabular}
\end{adjustbox}
\end{table}

% \subsection{A Subsection Sample}
% Please note that the first paragraph of a section or subsection is
% not indented. The first paragraph that follows a table, figure,
% equation etc. does not need an indent, either.

% Subsequent paragraphs, however, are indented.

% \subsubsection{Sample Heading (Third Level)} Only two levels of
% headings should be numbered. Lower level headings remain unnumbered;
% they are formatted as run-in headings.

% \paragraph{Sample Heading (Fourth Level)}
% The contribution should contain no more than four levels of
% headings. Table~\ref{tab1} gives a summary of all heading levels.

% \noindent Displayed equations are centered and set on a separate
% line.
% \begin{equation}
% x + y = z
% \end{equation}
% Please try to avoid rasterized images for line-art diagrams and
% schemas. Whenever possible, use vector graphics instead (see
% Fig.~\ref{fig1}).

% \begin{figure}
% \includegraphics[width=\textwidth]{fig1.eps}
% \caption{A figure caption is always placed below the illustration.
% Please note that short captions are centered, while long ones are
% justified by the macro package automatically.} \label{fig1}
% \end{figure}

% %
% \bibliographystyle{splncs04}
% \bibliography{references}
\pagebreak
\bibliographystyle{splncs04}
\bibliography{mybibliography}

\begin{thebibliography}{10}
\providecommand{\url}[1]{\texttt{#1}}
\providecommand{\urlprefix}{URL }
\providecommand{\doi}[1]{https://doi.org/#1}

\bibitem{amiot2009new}
Amiot, E.: New perspectives on rhythmic canons and the spectral conjecture.
  Journal of Mathematics and Music  \textbf{3}(2),  71--84 (2009)

\bibitem{amiot2011structures}
Amiot, E.: Structures, algorithms, and algebraic tools for rhythmic canons.
  Perspectives of new music  \textbf{49}(2),  93--142 (2011)

\bibitem{anders2018compositions}
Anders, T.: Compositions created with constraint programming. The Oxford
  Handbook of Algorithmic Music  (2018)

\bibitem{anders2011constraint}
Anders, T., Miranda, E.R.: Constraint programming systems for modeling music
  theories and composition. ACM Computing Surveys (CSUR)  \textbf{43}(4),
  1--38 (2011)

\bibitem{auricchio2021integer}
Auricchio, G., Ferrarini, L., Lanzarotto, G.: An integer linear programming
  model for tilings. arXiv preprint arXiv:2107.04108  (2021)

\bibitem{bailleux2003efficient}
Bailleux, O., Boufkhad, Y.: Efficient cnf encoding of boolean cardinality
  constraints. In: International conference on principles and practice of
  constraint programming. pp. 108--122. Springer (2003)

\bibitem{chemillier2001two}
Chemillier, M., Truchet, C.: Two musical csps. In: CP 01 Workshop on Musical
  Constraints. pp.~1--1 (2001)

\bibitem{courtot1990constraint}
Courtot, F.: A constraint-based logic program for generating polyphonies. Ann
  Arbor, MI: Michigan Publishing, University of Michigan Library (1990)

\bibitem{ebciouglu1988expert}
Ebcio{\u{g}}lu, K.: An expert system for harmonizing four-part chorales.
  Computer Music Journal  \textbf{12}(3),  43--51 (1988)

\bibitem{fripertinger2005remarks}
Fripertinger, H., Reich, L., et~al.: Remarks on rhythmical canons. Citeseer
  (2005)

\bibitem{hiller1957musical}
Hiller~Jr, L.A., Isaacson, L.M.: Musical composition with a high speed digital
  computer. In: Audio Engineering Society Convention 9. Audio Engineering
  Society (1957)

\bibitem{imms-sat18}
Ignatiev, A., Morgado, A., Marques{-}Silva, J.: {PySAT:} {A} {Python} toolkit
  for prototyping with {SAT} oracles. In: SAT. pp. 428--437 (2018)

\bibitem{Matolcsi}
Kolountzakis, M.N., Matolcsi, M.: Complex {H}adamard matrices and the spectral
  set conjecture. Collectanea Mathematica  \textbf{Extra},  281–291 (2006)

\bibitem{kolountzakis2009algorithms}
Kolountzakis, M.N., Matolcsi, M.: Algorithms for translational tiling. Journal
  of Mathematics and Music  \textbf{3}(2),  85--97 (2009)

\bibitem{phdthesis}
Liang, J.H.: Machine Learning for SAT Solvers. Ph.D. thesis, University of
  Waterloo (Dec 2018)

\bibitem{pachet2001musical}
Pachet, F., Roy, P.: Musical harmonization with constraints: A survey.
  Constraints  \textbf{6}(1),  7--19 (2001)

\bibitem{philipp2015pblib}
Philipp, T., Steinke, P.: Pblib--a library for encoding pseudo-boolean
  constraints into cnf. In: International Conference on Theory and Applications
  of Satisfiability Testing. pp. 9--16. Springer (2015)

\bibitem{radicioni2007constraint}
Radicioni, D.P., Lombardo, V.: A constraint-based approach for annotating music
  scores with gestural information. Constraints  \textbf{12}(4),  405--428
  (2007)

\bibitem{Vuza}
Vuza, D.T.: Supplementary sets and regular complementary unending canons (part
  one, two, three, four). Perspectives of New Music  (1991-93)

\bibitem{wiggins1989representing}
Wiggins, G., Harris, M., Smaill, A.: Representing music for analysis and
  composition. In: Proceedings of the 2nd IJCAI AI/Music Workshop. pp. 63--71
  (1989)

\bibitem{zimmermann2001modelling}
Zimmermann, D.: Modelling musical structures. Constraints  \textbf{6}(1),
  53--83 (2001)

\end{thebibliography}

\end{document}